\documentclass[12pt]{amsart}
\usepackage{amssymb,amsmath}

\numberwithin{equation}{section}
\def\T{\text}

\def\1#1{\overline{#1}}
\def\2#1{\widetilde{#1}}
\def\3#1{\widehat{#1}}
\def\4#1{\mathbb{#1}}

\def\5#1{\frak{#1}}
\def\6#1{{\mathcal{#1}}}
\def\C{{\4C}}
\def\R{{\4R}}

\def\Z{{\4Z}}

\begin{document}
\abstract
We give positive answer to a conjecture by Agranovsky. A continuous function on the sphere which has separate holomorphic extension along the complex lines which pass through three non aligned interior points, is the trace of a holomorphic function in the ball.
\newline
MSC: 32F10, 32F20, 32N15, 32T25 
\endabstract
\title[Separate holomorphic extension...]{Separate holomorphic extension along lines and holomorphic extension from the sphere to the ball} 
\author[L.~Baracco]{Luca Baracco}
\address{Dipartimento di Matematica, Universit\`a di Padova, via 
Trieste 63, 35121 Padova, Italy}
\email{baracco@math.unipd.it}
\maketitle
\def\Giialpha{\mathcal G^{i,i\alpha}}
\def\cn{{\C^n}}
\def\cnn{{\C^{n'}}}
\def\ocn{\2{\C^n}}
\def\ocnn{\2{\C^{n'}}}
\def\const{{\rm const}}
\def\rk{{\rm rank\,}}
\def\id{{\sf id}}
\def\aut{{\sf aut}}
\def\Aut{{\sf Aut}}
\def\CR{{\rm CR}}
\def\GL{{\sf GL}}
\def\Re{{\sf Re}\,}
\def\Im{{\sf Im}\,}
\def\codim{{\rm codim}}
\def\crd{\dim_{{\rm CR}}}
\def\crc{{\rm codim_{CR}}}
\def\phi{\varphi}
\def\eps{\varepsilon}
\def\d{\partial}
\def\a{\alpha}
\def\b{\beta}
\def\g{\gamma}
\def\G{\Gamma}
\def\D{\Delta}
\def\Om{\Omega}
\def\k{\kappa}
\def\l{\lambda}
\def\L{\Lambda}
\def\z{{\bar z}}
\def\w{{\bar w}}
\def\Z{{\1Z}}
\def\t{{\tau}}
\def\th{\theta}
\emergencystretch15pt
\frenchspacing
\newtheorem{Thm}{Theorem}[section]
\newtheorem{Cor}[Thm]{Corollary}
\newtheorem{Pro}[Thm]{Proposition}
\newtheorem{Lem}[Thm]{Lemma}
\theoremstyle{definition}\newtheorem{Def}[Thm]{Definition}
\theoremstyle{remark}
\newtheorem{Rem}[Thm]{Remark}
\newtheorem{Exa}[Thm]{Example}
\newtheorem{Exs}[Thm]{Examples}
\def\Label#1{\label{#1}}
\def\bl{\begin{Lem}}
\def\el{\end{Lem}}
\def\bp{\begin{Pro}}
\def\ep{\end{Pro}}
\def\bt{\begin{Thm}}
\def\et{\end{Thm}}
\def\bc{\begin{Cor}}
\def\ec{\end{Cor}}
\def\bd{\begin{Def}}
\def\ed{\end{Def}}
\def\br{\begin{Rem}}
\def\er{\end{Rem}}
\def\be{\begin{Exa}}
\def\ee{\end{Exa}}
\def\bpf{\begin{proof}}
\def\epf{\end{proof}}
\def\ben{\begin{enumerate}}
\def\een{\end{enumerate}}
\def\dotgamma{\Gamma}
\def\dothatgamma{ {\hat\Gamma}}

\def\simto{\overset\sim\to\to}
\def\1alpha{[\frac1\alpha]}
\def\T{\text}
\def\R{{\Bbb R}}
\def\I{{\Bbb I}}
\def\C{{\Bbb C}}
\def\Z{{\Bbb Z}}
\def\B{{\Bbb B}}
\def\P{{\Bbb P}}
\def\S{{\Bbb S}}
\def\Fialpha{{\mathcal F^{i,\alpha}}}
\def\Fiialpha{{\mathcal F^{i,i\alpha}}}
\def\Figamma{{\mathcal F^{i,\gamma}}}
\def\Real{\Re}
%
%
%
\section{  Introduction}
The problem of describing families  of discs which suffice for testing analytic extension of a function $f$ from the sphere $\partial\B^2$ to the ball $\B^2$ has a long history. For $f$ continuous on $\partial \B^2$, Agranovsky-Valski \cite{AV71} use all the lines, Agranovki-Semenov \cite{AS91} the lines through an open subset $D'\subset\B^2$, Rudin \cite{R} the lines tangent to a concentric subsphere $ B^2_{\frac12}$, Baracco--Tumanov-Zampieri the lines tangent to any strictly convex subset $D'\subset\subset \B^2$. 
There are many other contributions such as \cite{A07},\cite{S77}, \cite{GS91} just to mention a few.
It is a challenging attempt to reduce the number of parameters in the testing families. However, one encounters an immediate constraint: lines which meet a single point $z_o\in \B^2$ do not suffice. Instead, two interior points or a single boundary point suffice: Agranovsky \cite{A09} and Baracco \cite{B09}.
However, in these last two results, the reduction of the testing families is compensated by an assumption of extra initial regularity: $f$ is assumed to be real analytic.
Globevnik \cite{G09} shows that, for  two points,  $C^\infty$-regularity still suffices, but  $C^k$ does not. This suggests that  holomorphic extension is a good balance between reduction of testing families and improvement of  initial regularity. And in fact, it is showed  here, that for $f\in C^0$ three not on the same line points suffice. 
Here is our result.
\bt 
\Label{t1.1}
Let $f$ be a continuous function on the sphere $\partial \B^2$ which extends holomorphically along any complex line in $\B^2$ which encounters the set consisting of $3$ points not on the same line. Then, $f$ extends holomorphically to $\B^2$.
\et 
The proof follows in Section 2 below.
It shows that, the result should hold for a ball of general dimension $\B^n$. In this case, $n+1$ points in generic position should suffice.
We first introduce some terminology. A straight disc $A$ is the intersection of a straight complex line with $\B^2$;  $\Bbb PT^*\C^2$ is the cotangent bundle with projectivized fibers, and $\pi$ the projection on the base space; $\Bbb P T^*_{\partial \B^2}\C^2$ the projectivized conormal bundle to $\partial \B^2$ in $\C^2$. 
It is readily seen that the straight discs $A$ of the ball are the geodesics of the Kobayashi metric, or, equivalently,  the so called ``stationary discs" (cf. Lempert \cite{L81}). These are the discs endowed   with a meromorphic lift $A^*\subset \Bbb PT^*\C^2$ with a simple pole attached to $ T^*_{\partial \B^2}\C^2$, that is, satisfying $\partial A^*\subset \Bbb P T^*_{\partial \B^2}\C^2$.
We fix three points $P_j$, $j=1,2,3$ in $\Bbb B^2$ and  consider a set, indexed by $j$, of $(2)$-parameter families of  straight discs $A^j$ passing through $P_j$. We define $M_j$ to be the union of the lifts of the family with index $j$.
 The set $M_j$ is  generically a CR manifold  with CR dimension 1 except at the points  that project over $P_j$; we denote by $M_j^{\T{reg}}$ the complement of this set. The boundary of $M_j$ concides with $\Bbb P T^*_{\partial\Bbb B^2}\C^2$ which is maximal totally real in $\Bbb PT^*\C^2$.
Here is the central point of our construction. Though the function $f$, in the beginning of the proof,  is not extendible
to $\B^2$ as a result of the separate extensions to the $A$'s, nevertheless it is
 naturally lifted to a function $F$ on $M_j$ by gluying the bunch of separate holomorphic extensions to the lifts $A^*$'s.
 This is defined by
\begin{equation*}
F(z,[\zeta])=f_{A_{(z,[\zeta])}}(z),
\end{equation*}
where $A_{(z,[\zeta])}$ is the unique stationary disc whose lift $A^*_{(z,[\zeta])}$ passes through $(z,[\zeta])$. The  crucial point here is that the $A$'s may overlap on $\C^2$ but the $A^*$'s do not in $\mathbb{P} T^*\C^2$. The function $F$ is therefore well defined and  CR on $M_j^{\T{reg}}$.


\section{Proof of Theorem~\ref{t1.1}}
The proof consists of several steps.
We start by collecting  some easy computations. We identify $\P T^*\C^2\simeq \C^2\times\C\P_1\simeq \C^3$ with coordinates $(z_1,z_2)\in\C^2$ and $z_3= \frac{\zeta_2}{\zeta1}\in\C\P_1$. Let $M_0$ be the collection of the lifts of the discs through $0$.
\bl 
Let $A^*_0$ be the (unique) disc of $M_0$ which projects over the $z_1$-axis. Then, $A^*_0$, identified to a disc of $\C^3$, has two holomorphic lifts to $T^*\C^3$ attached to $T^*_{M_0}\C^3$. Their components are parametrized by $z_1\mapsto (0,-\frac 1{z_1},1)$ and $z_1\mapsto (0,\frac 1{iz_1},\frac 1i)$respectively.
\el
\bpf 
First, we notice that for any $z=(z_1,z_2)\in \B^2$ the disc $\tau\mapsto \tau\frac z{\|z\|}$ is the only passing through $z$ and $0$. The lift attached to the projectivized conormal bundle of this disc  is  the constant $[\bar{z}]$. We have
$$M_0=\{ (z;[\bar{z}])\  z\in\B^2\setminus 0\}\cup\{ (0;[\zeta])\  \forall [\zeta]\in \C\P_1\}. $$
Clearly $M_0$ (or more precisely $M_0^{\text{reg}}$) has equation $r:z_3-\frac{\bar{z}_2}{\bar{z}_1}=0$. 
In particular the lift of $A_0$  to $\Bbb P T^*\C^2$ is $A_0^*(\tau)=((\tau,0); [1,0])$ which in coordinates is expressed by $A^*_0(\tau)=(\tau,0,0)$. Since $M_0$ is Levi flat, the space of holomorphic lifts contained in $T^*M_0$  has dimension two. For instance a basis for the space of lifts is given by
\begin{equation}
\omega_1(z_1,z_2)=\partial\Re r= \left(\frac{z_2}{z^2_1},-\frac 1 {z_1},1\right) \text{and } \omega_2(z_1,z_2)=\partial\Im r=\frac 1i\left(-\frac{z_2}{z^2_1},\frac 1 {z_1},1\right).
\end{equation}

In particular, along $A_0^*$ the conormal bundle to $M_0$ is generated by $\omega_1(z_1,0)=(0,\frac {-1}{z_1},1)$ and $\omega_2(z_1,0)=(0,\frac 1{iz_1},\frac 1i)$. As one can readily note both sections are holomorphic along $A^*_0$ and they are exactly the lifts of $A^*_0$ to the conormal bundle of $T^*_{M_0}\C^3$. 

\epf

\br Note that if in the  lemma above we consider the union of the lifts of discs passing through the point $P_{\zeta_0}=(\zeta_0,0)$ the manifold resulting $M_{\zeta_0}$ still contains $A_0^*$ (i.e. the $z_1$ axis) and along the boundary of $A_0^*$ we have 
$TM_0|_{\partial A^*_0}=TM_{\zeta_0}|_{\partial A^*_0}$ and thus also $T^*_{M_0}\C^3|_{\partial A^*_0}=T^*_{M_{\zeta_0}}\C^3|_{\partial A^*_0}$.
From this equality we have that if $\tilde{\omega}_1$, $\tilde{\omega}_2$ is a basis of lifts of $A^*_0$ to the conormal bundle to $M_{\zeta_0}$, then this is related to the basis $\omega_1,\omega_2$ by
\begin{equation}
\label{conormal}
\text{Span} \{\tilde{\omega}_1, \tilde{\omega}_2\}|_{\partial A^*_0} =\text{Span} \{\omega_1, \omega_2\}|_{\partial A^*_0}.
\end{equation}
Combination of (\ref{conormal}) with the fact that singularity of $\tilde{\omega_1}, \tilde{\omega_2}$ must now be located at $\zeta_0$ yields a choice of holomorphic basis as
$\tilde{\omega_1}(z_1)=\left(0,-\frac{1}{(z_1-\zeta_0)},\frac{1}{(1-z_1\bar\zeta_0)}\right)$ and $\tilde{\omega_2}(z_1)=\left(0,\frac 1{i(z_1-\zeta_0)},\frac{1}{i(1-z_1\bar\zeta_0)}\right)$.
\er
Before the proof of our main theorem we need a preliminary crucial result
\bt
Let $P_1,P_2\in\B^2$ be two distinct points inside the ball and let $f:\partial \B^2\rightarrow \C$ be a continuous function such that $f$ extends holomorphically along every complex straight line passing through either $P_1$ or $P_2$. Then for any such disc $A$, except the one passing through both points, the lifted function $F$ extends holomorphically to a neighborhood of $A^*\setminus \pi^{-1}(P_j)$ where $j$ is $1$ or $2$ according to $P_1\in A$ or  $P_2\in A$. 
\et
\bpf
It is not restrictive to assume that the disc $A$ is the $z_1$ axis, that $P_2=(0,z_2)$ and that $P_1=(\zeta_0,0)$.
We note that $M_1$ and $M_2$ intersect transversally along the boundary of $A^*$.
Let $P=(\zeta,0)$ be a point of the boundary of $A$ and $P^*=(\zeta,0,0)$ be the corresponding point on $A^*$. $P^*$ lies in the common boundary of $M_1$ and $M_2$. Let $v_\zeta $ be a tangent vector to $T_{P^*}M_2\setminus T_{P^*}E$ which points inside $M_2$. The equivalence class $[v_\zeta]$ in the vector spaces quotient $\frac{T_{P^*}\C^3}{ T_{P^*}M_1}$ is called the pointing direction of $M_2$ with respect to $M_1$. We say in this case that $F$ extends at $P^*$ in direction $[v_\zeta]$. Let $Q^*=(\zeta_Q,0,0)$ be a point of $A^*$ ($\zeta_Q\neq\zeta_0$). Following \cite{T97} by effect of the extension of $F$ at $P^*$ in direction $[v_\zeta]$ we have extension at $Q^*$ in direction $[w_\zeta]\in \frac{ T_{Q^*}\C^3}{T_{Q^*}M_1}$. The relation of $[w_\zeta]$ with the initial $[v_\zeta]$ is expressed by means of contraction with the holomorphic basis of lifts $\tilde{\omega_1},\tilde{\omega_2}$:
\begin{equation}
\label{c0}
 \Re  \langle\tilde{\omega}_1(\zeta),v_\zeta\rangle =\Re  \langle \tilde{\omega}_1(\zeta_Q),w_\zeta\rangle \text{ and }\Re  \langle\tilde{\omega}_2(\zeta),v_\zeta\rangle =\Re  \langle \tilde{\omega}_2(\zeta_Q),w_\zeta\rangle .
\end{equation}
In other words the directions of $CR$ extendibility, which are  vectors in the normal bundle $\frac{T\C^3}{TM_1}$, are constant in the system dual to $\{\tilde{\omega}_1,\tilde{\omega}_2\}$.  

  We first compute the pointing direction of $M_2$ at the point $P^*$. To this end we first compute the disc passing through $P_2$ and $P$  which is
$$A_{P_2P}(\tau)=(\frac{|z_2|^2\zeta}{1+|z_2|^2},\frac{z_2}{1+|z_2|^2})+\frac{\tau}{1+|z_2|^2}(\zeta,-z_2);$$
note that $A_{P_2P}(1)=P$. The lift component of $A_{P_2P}$ is
$$A^*_{P_2P}=[|z_2|^2\bar\zeta\tau+\bar\zeta,\bar z_2\tau-\bar z_2],$$
and dividing the second component by the first we get that the $A^*_{PP_2}$'s coordinates in $\C^3$ are
$$  \left( (\frac{|z_2|^2\zeta}{1+|z_2|^2}+\frac{\tau}{1+|z_2|^2}\zeta,\frac{z_2}{1+|z_2|^2}-\frac{\tau z_2}{1+|z_2|^2},\frac{\bar z_2(\tau-1)}{\bar\zeta(|z_2|^2\tau +1)}\right).$$
The pointing direction of $M_2$ at $P$ is
$$v_{\zeta}=-\partial_\tau A^*_{P_2P}(1)=\frac{-1}{1+|z_2|^2}(\zeta,-z_2,\frac{\bar z_2}{\bar\zeta}).$$
We have 
\begin{equation}
\label{c1}
 \Re  \langle\tilde{\omega}_1(\zeta),v_\zeta\rangle=\frac{-1}{1+|z_2|^2} \Re  \frac{z_2}{\zeta-\zeta_0}
\end{equation}
and
\begin{equation}
\label{c2}
 \Re \langle\tilde{\omega}_2(\zeta),v_\zeta\rangle=\frac{-1}{1+|z_2|^2} \Im  \frac{z_2}{\zeta-\zeta_0}.
\end{equation}
The first members of (\ref{c1}) and (\ref{c2}) express the components in the normal bundle to $M_1$ of $w_\zeta$ with respect to the dual basis of $\omega_1(\zeta_Q),\omega_2(\zeta_Q) $.  
By letting $\zeta$ vary in $\partial A$ we see that $[w_\zeta]$ sweeps all the directions in $\frac{T\C^3}{TM_1}|_{Q^*}$. Therefore, by the edge of the wedge theorem, $F$ extends holomorphically to a neighborhood of $Q^*$ and, by propagation, to a neighborhood of any other point of $A^*$ except from the point over $P_1$ where the $CR$ singularity is located.

\epf    
We are ready for the proof of Theorem~\ref{t1.1}

\noindent{\it End of Proof of Theorem~\ref{t1.1}} 
 Let $A_0$ be the disc passing through $P_1$ and $P_3$. Then in particular $P_2\notin A_0$. Applying the theorem  above we get that $F$ extends holomorphically to a neighborhood of $A_0^*\setminus\{P_1\}$. By repeating the same argument we see that $F$ extends to a neighborhood of $A_0^*\setminus\{P_3\}$. Therefore $F$ extends to a full neighborhood of $A^*_0$. For any other straight line $A$ through $P_1$ we can say that $F$ extends holomorphically to a neighborhood of $A^*\setminus{P_1}$. By applying the continuity principle to the family of discs formed by $A^*_0$ and all the discs through $P_1$, we get that $F$ extends holomorphically to a set of the form $V\times \C\P^1_\C$ where $V$ is a neighborhood of $P_1$. Therefore $F$ does not depend on the second argument and it is therefore naturally defined on the projection of the collection of all the $A^*$'s, that is, on $\B^2$. 

\hskip11cm$\Box$

\end{document}